\documentclass{article}
\usepackage{amssymb}
\usepackage{amsmath}
\newtheorem{theorem}{Theorem}
\newtheorem{lemma}{Lemma}\newtheorem{proposition}{Proposition}\newtheorem{definition}{Definition}
\newtheorem{corollary}{Corollary}
\newcommand{\La}{ \Lambda}
\newcommand{\R}{{\mathbb R}}

\newcommand{\N}{{\mathbb N}}

\newcommand{\C}{{\mathbb C}}

\usepackage{lipsum}

\newcommand\blfootnote[1]{%
  \begingroup
  \renewcommand\thefootnote{}\footnote{#1}%
  \addtocounter{footnote}{-1}%
  \endgroup
}

\begin{document}
\begin{large}

\title{On Cartwright's  theorem}

\author{Natalia Blank  and  Alexander Ulanovskii}
\date{}\maketitle

\begin{abstract}
We present a characterization of  sets  for which Cartwright's theorem holds true.  The connection is discussed between these sets and sampling sets for  entire functions of exponential type.
\end{abstract}


\section{Cartwright's and Beurling's Theorems}

\subsection{Cartwright's Theorem} An entire function $f(z)$ is said to be of  exponential type if there exist positive numbers $c$ and $C$
such that
\begin{equation}\label{0}
|f(z)|<Ce^{c|z|},\ \ \mbox{ for all $z\in \C$}.\end{equation}
The type of $f(z)$ is defined as  the infimum over all $c$ which can be used. The function $f$ is of exponential type zero, if (\ref{0}) holds for every $c>0$.

We will denote by  $E_\sigma$  the class of all entire functions of exponential type $\leq \sigma,$ and by $E_{<\sigma}$ the class of all entire functions of exponential type $< \sigma.$

  A  theorem of M. L. Cartwright states:

 \begin{theorem}[\cite{c}]  If a function $f\in E_{<\pi}$ is bounded on the set of integers then $f$ is bounded on the real axis.
\end{theorem}

The example $f(z)=z\sin(\pi z)$ shows that Theorem 1 ceases to be true  for the functions of exponential type $\pi$.

A  set $\La\subset\R$ is called uniformly discrete (u.d.) if \begin{equation}\label{ud} d(\La):=\inf_{\lambda,\lambda'\in\La,\lambda\ne\lambda'}|\lambda-\lambda'|>0. \end{equation} The number $d(\La)$ is called the separation constant of $\La.$ Denote by $D^-(\La)$ the lower uniform density of $\La$: $$ D^-(\La):=\lim_{r\to\infty}\min_{x\in\R}\frac{\#(\La\cap(x-r,x+r))}{2r}. $$

%\medskip\noindent{\bf Definition 1}.

\begin{definition}  We say that $\La\subset\R$ is a Cartwright  set (CS) for a class $M$ of entire functions, if there is no function   $f\in M$ which is bounded on $\La$ and unbounded
on the real axis.\end{definition}

  Theorem 1 states that the set of integers is a CS for the class $E_{<\pi}$.

We prove

 \begin{theorem} A set $\La\subset\R$ is a CS for $E_{<\sigma}$ if and only if it contains a u.d. subset $\La^*$ satisfying $D^-(\La^*)\geq \sigma/\pi.$\end{theorem}

 Since the subject is classical nowadays, this result may be known to the experts. However, to the best of our knowledge, it has not been recorded in the literature.

 Cartwright's theorem has inspired a number of different results. In particular, it is shown  in  \cite{b} and   \cite{df} (see also literature therein),  that   certain ``perturbations" of the integers are CSs for $E_{<\pi}$.

Papers   \cite{a}, \cite{b1}, \cite{l1} and \cite{m} adopt a more general approach: The condition of boundedness of $f$ on $\La$ is dropped, and   results are obtained concerning the growth of the function from its growth on  $\La$.

 See  \cite{bs}, \cite{lm},  \cite{ou} (and the literature therein) for some estimates of   $\|f\|_\infty$ for functions $f\in E_\sigma$ satisfying $\|f|_\La\|_\infty\leq 1$. Here
$$
\|f\|_\infty:=\sup_{x\in\R}|f(x)|, \ \|f|_\La\|_\infty:=\sup_{\lambda\in\La}|f(\lambda)|.
$$

  See also \cite{lf}, \cite{lo} and \cite{ln} for  some  multi-dimensional Cartwright-type results.

 There have been a number of other results related to Theorem 1, which we do not mention here.

 It seems plausible that  Theorem 2 can be deduced from  the results of \cite{a} or \cite{m}. In our proof we use A. Beurling's results  from \cite{beu}.

\medskip

\subsection{Beurling's Sampling Theorem}   Let $B_\sigma$ denote the subclass of $E_\sigma $ of functions $f$ bounded on the real axis. It is well--known that every function  $f\in B_\sigma$ satisfies the inequality (see, for example \cite{levin}, Lec. 6, Theorem 3)
\begin{equation}\label{ber}
|f(x+iy)|\leq\|f\|_\infty e^{\sigma|y|}, \quad \mbox{ for all } x,y\in\R.
\end{equation}The spaces $B_\sigma$ are usually called Bernstein spaces.

A set $\La$ is called a sampling set (SS) for $B_\sigma$, if there is a constant $K$ such that
\begin{equation}\label{s}
\|f\|_\infty\leq K\|f|_\La\|_\infty,\quad \mbox{ for all } f\in B_\sigma.
\end{equation}
Denote by $K(\La, B_\sigma)$ the infimal  $K$ in (\ref{s}), and set $K(\La,B_\sigma)=\infty$ if $\La$ is not an SS for $B_\sigma$.
 The constant $K(\La,B_\sigma)$ is  called the sampling bound for   $B_\sigma$.

A classical result of   Beurling  describes the SSs for $B_\sigma $:

 \begin{theorem}[\cite{beu}]  A  set $\La\subset\R$ is an SS for $B_\sigma$ if and only if it contains a u.d. subset $\La^*$ satisfying  $D^-(\La^*)>\sigma/\pi.$\end{theorem}

\section{Sampling Sets and  Cartwright Sets for $E_\sigma$}

Let us extend the notion of sampling set.

\begin{definition} We say that a set $\La\subset\R$ is an SS for a class $M$ of entire functions,
   if there is a constant
$K$ such that for every $f\in M$ the inequality
\begin{equation}\label{3}
 \max_{|x|\leq r}|f(x)|\leq K \max_{\lambda\in\La,|\lambda|\leq r}|f(\lambda)|
\end{equation}holds on an unbounded set of $r\in(0,\infty)$ (this set depends on $f$).
\end{definition}

This definition  means that, up to a multiplicative constant, no function $f\in M$ can grow faster along the  real axis than it grows  along $\La$.
When $M=B_\sigma,$ this definition  coincides with the classical one.

 In what follows, we denote by $K(\La,M)$ the infimal  $K$ in (\ref{3}), and we set $K(\La,M)=\infty$ if $\La$ is not an SS for $M$.

 Since $B_\sigma\subset E_\sigma$, we have   $K(\La, B_\sigma)\leq K(\La,E_\sigma).$
In fact, these two constants are equal:

 \begin{theorem}  A set $\La\subset\R$ is an SS for $E_\sigma$ if and only if it is an SS for $B_\sigma$. Moreover, $K(\La,E_\sigma)=K(\La, B_\sigma)$.\end{theorem}

 It is clear from the definitions above that every SS for   a class $M$  is also a CS for $M$. We show that for the classes $E_\sigma$ the opposite is also true:

 \begin{theorem} A  set $\La\subset\R$ is  a CS for $E_\sigma$ if and only if it is  an SS for $E_\sigma.$ \end{theorem}

 Theorems 3,4 and 5 yield

  \begin{corollary}
  A  set $\La\subset\R$ is  a CS for $E_\sigma$ if and only if it contains a u.d. subset $\La^*$ satisfying  $D^-(\La^*)>\sigma/\pi.$
  \end{corollary}

 Theorem 2 above is an immediate consequence of this corollary.

\medskip\noindent{\bf Remark 1}. It might be interesting to compare Theorem 2 with the result in \cite{pm} which gives a characterization of the  P\'{o}lya sets. A set $\La\subset\R$ is called a P\'{o}lya set if the constants are the only entire functions of exponential type zero which are bounded on $\La$. Since there is no non-constant function of zero type  bounded on $\R$,  the P\'{o}lya sets are therefore  the  Cartwright  sets for the class of  entire functions of exponential type zero.  However, the characterization of CSs for functions of finite exponential type involves the lower uniform density, while the characterization of P\'{o}lya sets involves the so-called inner Beurling-Malliavin density, see \cite{pm}.

\medskip\noindent{\bf Remark 2}.
One may wish to describe CSs for other classes of entire functions. Here we present a result in this direction.

 Given a natural number $n$ and a positive number $\sigma$, we denote by $E_{<\sigma,n}$ the space of all entire functions $f$ satisfying the inequality
$$
|f(z)|\leq Ce^{c|z|^n}, \ \ z\in\C,
$$with some $0<c<\sigma$ and $C>0$. Clearly, $E_{<\sigma,1}=E_{<\sigma}$.

The following result extends Theorem 2 to the classes $E_{<\sigma,n}$:

 \begin{theorem}  A set $\Gamma\subset\R$ is a CS for $E_{<\sigma,n}$ if and only if $\Gamma^n$ is  a CS for $E_{<\sigma}$, where $$
\Gamma^n:=\{\gamma^n: \gamma\in\Gamma, \gamma\geq0\}\cup\{-|\gamma|^n: \gamma\in\Gamma, \gamma<0\}.
$$ \end{theorem}

By Theorem 2, it follows  that $\Gamma$ is a CS for $E_{<\sigma,n}$ if and only if the set $\Gamma^n$ contains a u.d. subset $\La$ satisfying $D^-(\La)\geq\sigma/\pi.$

\section{Auxiliary Result}

\noindent{\bf 1}. We will need several results from  \cite{beu}.

Given two sets $\La,\La^*$ such that $\La^*\subset\La$, it is clear that  $K(\La^*, B_\sigma)\geq K(\La, B_\sigma)$.  The next proposition shows that for every sampling set $\La$, one may pick up a uniformly discrete (u.d.) subset $\La^*$ without changing much the sampling bound:

\begin{proposition}[\cite{beu}] Suppose $\La$ is an SS for $B_\sigma.$ Then for every $\varepsilon>0$ there is a u.d. subset $\La^*\subset\La$ such that
 $K(\La^*, B_\sigma)<K(\La, B_\sigma)+\varepsilon.$\end{proposition}

 For a given closed set $\La$, let $\La(t)=\La+[-t,t]$ denote the set of  points with distance $\leq t$ from $\La$. The Hausdorff %Fr\'{e}chet
 distance between two closed sets $\La$ and $\Gamma$ is the smallest number $t$ so that $$
\La\subset \Gamma(t), \Gamma\subset\La(t).
$$We say that a sequence of closed sets  $\La_n$ converges weakly to a closed set $\La$, if for every closed interval $I=[a,b], a,b\not\in\La$, the distance between $\La_n\cap I$ and $\La\cap I$ tends to zero as $n\to\infty.$

The next proposition plays a key role in Beurling's theory:

\begin{proposition}[\cite{beu}] Assume u.d. sets $\La_j$ converge weakly to a u.d. set $\La$. Then $$K(\La,B_\sigma)\leq \lim\inf_{j\to\infty}K(\La_j,B_\sigma).$$ \end{proposition}

The following elementary lemma is stated without proof:

\begin{lemma} (i) Every sequence of  sets  $\La_j$ satisfying $d(\La_j) > d >0, j\in\N$, where $d(\La_j)$ is the separation constant in (\ref{ud}), has a subsequence  converging weakly to some set $\La$ satisfying $d(\La) > d$.

(ii)  The translations of a set $\La$ do not change the sampling bound, i.e. for every real number $a$ we have
$K(\La-a,B_\sigma)=K(\La,B_\sigma)$.
\end{lemma}

\noindent{\bf 2}.  We will also need  a variant of Theorem~4 for functions analytic in the  right half-plane $\C_r:=\{z: \Re z> 0\}$:

\begin{theorem}
Let $C>0,\sigma>0$ and let   $\La\subset\R$ be an SS for  $B_\sigma.$ Assume a function $f$ is  analytic  in the half-plane $\C_r$, continuous in its closure and satisfies
\begin{equation}\label{f}
|f(x+iy)|<Ce^{\sigma\sqrt{x^2+ y^2}}, \ \ \mbox{ for all } x\geq0, y\in\R,
\end{equation}
and \begin{equation}\label{ff}\limsup_{x\to+\infty}|f(x)|=\infty.\end{equation}
Then for every $K>K(\La, B_\sigma)$ the inequality
\begin{equation}
\max_{0\leq x\leq r}|f(x)|\leq K\max_{\lambda\in\La,0\leq \lambda\leq r}|f(\lambda)|
\end{equation}holds on an unbounded set of $r\in(0,\infty)$.
\end{theorem}

\noindent{\bf Proof of Theorem 7} (By contradiction).
  Suppose that the assumptions of Theorem 7 are fulfilled and that (8)  is not true: There are numbers $K>K(\La,B_\sigma)$ and $r_0>0$ such that
  \begin{equation}\label{c}
  \max_{0\leq x\leq r}|f(x)|>K\max_{\lambda\in\La,0\leq \lambda\leq r}|f(\lambda)|, \ \ \mbox{ for all } r>r_0.
  \end{equation}

 By Proposition 1, there is a u.d. set $\La^*\subset\La$ satisfying $K(\La^*,B_\sigma)<K.$ So, in the rest of the proof we may assume that $\La$ is a u.d. set.

   Let us consider the auxiliary functions
 \begin{equation}\label{au}
  f_n(z):=e^{-\frac{1}{n}z\log z}f(z), \  \ \Re z>0, \ \ n\in\N.
   \end{equation}
 Similar functions were introduced by  Beurling (see the proof of Theorem 5 in \cite{beu}). We will use their  properties to show that (\ref{c}) leads to a contradiction.

  From (\ref{f}) and (\ref{au}), it follows that each $f_n$ tends to zero as $x\to+\infty$. Therefore, by (\ref{ff}), for every large enough  $n$, the function $|f_n(x)|$ attains its maximum on $(0,\infty)$. Let us denote by $x_n>0$  the smallest point  where $|f_n(x_n)|=\sup_{x>0}|f_n(x)|.$   Clearly, $x_n\to+\infty$ and $|f_n(x_n)|\to\infty$ as $n\to\infty.$

Observe that $x_n>r_0$ for all large enough $n$, where $r_0$ is the number in (\ref{c}). It follows from (\ref{c}) that for every $\lambda\in\La\cap(0,\infty)$ we have
  $$
  |f_n(\lambda)|=|f(\lambda)|e^{-\frac{1}{n}\lambda\log\lambda}<K\max_{0\leq x\leq\lambda}|f(x)|e^{-\frac{1}{n}\lambda\log\lambda}\leq$$
  $$K\max_{0\leq x\leq\lambda}|f_n(x)|\leq K|f_n(x_n)|.$$
Hence,
  \begin{equation}\label{cc}
|f_n(x_n)|>K\sup_{\lambda\in\La,\lambda>0}|f_n(\lambda)|.
  \end{equation}

 Using (\ref{f}) and (\ref{au}), one may easily check that  $f_n$ satisfies
 $$
|f_n(x+iy)|\leq Ce^{\sigma\sqrt{x^2+y^2}+\frac{\pi}{2n} |y|-\frac{1}{n}x\log x},   \mbox{ for all } x\geq0, y\in\R.
$$
Since $|f_n(x_n)|>C$ for all large enough $n$, this gives
$$
|f_n(iy)|\leq |f_n(x_n)|e^{(\sigma+\frac{\pi}{2n})|y|}, \  \mbox{ for all }  y\in\R.
$$
Recall that  $|f_n(x)|\leq |f_n(x_n)|, x>0.$ So, we may apply  the Phragmen-Lindel\"of  principle (see \cite{levin}, Lec. 6) to $f$ in the first and fourth quadrants to get the estimate
$$
|f_n(x+iy)|\leq |f_n(x_n)| e^{(\sigma+\frac{\pi}{2n}) |y|},\  \mbox{ for all } x\geq0, y\in\R.
$$

Put$$
h_{n}(z):=\frac{f_n(z+x_n)}{f_n(x_n)}, \ \ \Re z>-x_n.
$$Then
\begin{equation}\label{no}
h_{n}(0)=1, \ |h_{n}(x+iy)|\leq e^{(\sigma+\frac{\pi}{2n}) |y|},\  x>-x_n,y\in\R.
\end{equation}

Set $\Gamma_n:=\La-x_n$.
By (\ref{cc}),
$$
|h_{n}(\lambda)|<\frac{1}{K},  \ \ \lambda\in\Gamma_n, \lambda>-x_n.
$$By  Lemma 1 (i), there is a subsequence   $\Gamma_{n_j}$ converging weakly to some  $\Gamma$. Then by Lemma 1 (ii) and  Proposition 2, we have
\begin{equation}\label{o}
K(\Gamma,B_\sigma)\leq K(\La,B_\sigma)<K.
\end{equation}

Condition (\ref{no}) shows that for every compact $S$ in $\C$ there is an integer $n(S)$ such that  $\{h_{n}, n>n(S)\}$ is a normal family on $S$. Hence, it is clear that there is a subsequence $\{m_j\}\subset \{n_j\}$ such  that $h_{m_j}$ converge uniformly on  compacts to an entire function $h$ satisfying
$$
h(0)=1, \ |h(x+iy)|\leq e^{\sigma |y|},\ \ \mbox{ for all } x,y\in\R,
$$and $\| h|_\Gamma\|_\infty\leq 1/K.$ These inequalities show that  $h\in B_\sigma$ and  that $K(\Gamma, B_\sigma)\geq K,$ which  contradicts  (\ref{o}).
\section{Proofs of Theorems 4 and 6}

\subsection{Proof of Theorem 4} It suffices to check that $K(\La,E_\sigma)\leq K(\La,B_\sigma)$. In other words, it suffices to check that inequality (\ref{3}) holds on an unbounded set of $r>0$ for every function $f\in E_\sigma$ and every constant $K>K(\La,B_\sigma)$.

If $f$ is bounded on $\R$, then $f\in B_\sigma$. Hence, (\ref{3}) holds by the definition of sampling bound $K(\La,B_\sigma)$.

Assume $f$ is  unbounded on $(0,\infty)$ and bounded on $(-\infty,0)$. Then (\ref{3}) immediately follows from Theorem 7.

To establish (\ref{3}) in the case when $f$ is  unbounded on $(-\infty,0)$ and bounded on $(0,\infty),$ we apply Theorem 7 to the function $f(-z)$.

Now assume $f$ that is unbounded on both $(-\infty,0)$ and  $(0,\infty)$.
Let us assume that (\ref{3}) is not true, i.e. for some $K>K(\La, B_\sigma)$ and $r>r_0$ we have
\begin{equation}\label{kk} \max_{|x|\leq r}|f(x)|> K\max_{|\lambda|\leq r,\lambda\in\La}|f(\lambda)|,  \ \ \ r>r_0.\end{equation}

Consider the functions
$$
g_n(x):=|f(x)|e^{-\frac{1}{n}|x|\log|x|}, \ \ x\in\R.
$$As in the proof of Theorem 7, for all large enough $n$ we denote by  $x_n$ a point with the smallest absolute value such that
$|g_n(x_n)|=\|g_n\|_\infty$.
We may  assume that  an infinite number of $n$  is positive  (otherwise, we consider $g_n(-x)$). Then, for such $n$ we consider the functions  $f_n$ defined in the proof of Theorem 7. By (\ref{kk}), it is easy to check that $f_n$ satisfies (\ref{cc}). Then one may repeat the proof of Theorem 7 to show that (\ref{kk}) leads to a contradiction.

\subsection{Proof of Theorem 6}  (i) Let $\Gamma\subset\R$ be such that $\La:=\Gamma^n$ is a CS for $E_{<\sigma}$. Then, by Theorem 4, $\La$ is an SS for every space $B_s, s<\sigma.$

Let us show that  $\Gamma$ is a CS for $E_{<\sigma,n}$. It suffices to show that if a function $f\in E_{<\sigma,n} $ is unbounded on $\R$, then $f$ is also  unbounded on $\Gamma$.
We may assume that $f$ is unbounded on $(0,\infty)$. Set $g(z):=f(z^{1/n}), \Re z>0$. Using the definition of $E_{<\sigma,n}$, one may easily check that $g$ and $\La$ satisfy the assumptions of Theorem 7. Hence, $g$ is unbounded on $\La\cap(0,\infty)$, and so $f$ is unbounded on $\Gamma\cap(0,\infty)$.

If $f$ is unbounded on $(-\infty,0)$, we apply the argument above to $g(z):=f((-z)^{1/n}), \Re z<0$.

(ii) Assume $\La=\Gamma^n$ is not a CS for $E_{<\sigma}$. Then, by Theorem~5, there exists $g\in E_{<\sigma}$, which is bounded on $\La$ and unbounded on $\R$. Clearly, the function $f(z):=g(z^n)$ belongs to $E_{<\sigma,n}$,  $f$ is bounded on $\Gamma$  and unbounded on $\R$. We may assume that it is unbounded on $(0,\infty)$. Hence, $\Gamma$   is not a CS for $E_{<\sigma,n}$.

\section{Proof of Theorem 5}

Denote by $PW_\sigma$ the classical Paley--Wiener space of all Fourier transforms
$$
f(x)=\hat F(x):=\int_{-\sigma}^\sigma e^{-ixt}F(t)\,dt,
$$where $F\in L^2(-\sigma,\sigma).$
By the classical Paley--Wiener theorem (see \cite{levin}, Lec.10), $PW_\sigma$ consists of all entire functions of exponential type $\leq\sigma$ which are square-integrable on $\R$: $$PW_\sigma=L^2(\R)\cap E_\sigma.$$

We say that a set $\La$ is a uniqueness set for a function space $M$ if there is no non-trivial function $f\in M$ which vanishes on $\La$.

 We will need

\begin{lemma}
Suppose $\La\subset\R$ is not a uniqueness set for  $B_\sigma$. Then there exists $g\in E_\sigma$ which vanishes on $\La$ and is unbounded on $\R$.
 \end{lemma}

This lemma is a simple consequence of the following
\begin{proposition}[\cite{be}]
 Every incomplete systems of complex exponentials in $L^2(-\sigma,\sigma)$ is a subset of some complete and minimal system of exponentials.
 \end{proposition}

Denote by $E(\La)=\{e^{i\lambda t},\lambda\in\La\}$ the exponential system with frequencies in $\La$.

 Since
$\La$ is not a uniqueness set for $B_\sigma$,  there is a function $h\in B_\sigma$ which vanishes on $\La$.
 Then, for every $\lambda_0\in\La$, the function $h(z)/(z-\lambda_0)$ belongs to $PW_\sigma$. This means that
$E(\La\setminus\{\lambda_0\})$ is not complete in $L^2(-\sigma,\sigma)$. By Proposition~3, there exists $\Gamma\supset \La\setminus\{\lambda_0\}$ such that $E(\Gamma)$ is complete while every system $E(\Gamma\setminus\{\gamma_0\}),\gamma_0\in\Gamma,$ is not complete in $L^2(-\sigma,\sigma).$ 

We may assume that $\gamma_0\not\in \La\setminus\{\lambda_0\}$. Take  a function $F\in L^2(-\sigma,\sigma)$ which is orthogonal to $E(\Gamma\setminus\{\gamma_0\})$. Then its Fourier transform $f=\hat F$ belongs to $ PW_\sigma$, vanishes on $\La\setminus\{\lambda_0\}$ and   $(z-\gamma_0)f(z)\not\in PW_\sigma$.
Now, one may easily check that the function $g(z):=(z-\lambda_0)(z-\gamma_0)f(z)$ belongs to $E_\sigma$, vanishes on $\La$ and is unbounded on $\R$.

Observe that the proof of Proposition in \cite{be} is nonconstructive and based on the deep  de Branges theory of entire functions. At the end of this note we present a sketch of  a more elementary proof of Lemma 2.

\medskip\noindent{\bf Proof of Theorem 5}. Clearly, if $\Lambda$ is an SS for $E_\sigma$ then it is a CS for $E_\sigma$. So, it suffices to show that  if $\La$ is not an SS for $E_\sigma$, then  it is not a CS for $E_\sigma,$ i.e. there exists $f\in E_\sigma$ which is bounded on $\La$ and unbounded on $\R$.

We will consider three cases:

\medskip

\noindent{\bf 1}.  Assume $\La\subset\R$ is not a uniqueness set for  $B_\sigma$. Then Theorem 5 follows from Lemma 2.

\medskip
\noindent{\bf 2}.  Since $\La$ is not an SS for $B_\sigma$, there is a sequence of functions $f_n\in B_\sigma$ satisfying
\begin{equation}\label{sist}
\|f_n\|_\infty=1, \ \ \|f_n|_\La\|_\infty<\frac{1}{n^3}, \ \ n=1,2,...
\end{equation}

Assume that there is an interval $I=[-N,N], N>0,$ such that $\|f_n\cdot 1_N\|_\infty\not\to 0,$ as $n\to\infty.$ Here $1_N(x)$ is the indicator function of $[-N,N]$. In this case there is a point $a, |a|\leq N,$ and a subsequence
$n_j$ such that $|f_{n_j}(a)|>\delta>0.$ Since $f_n$ are uniformly bounded on $\R$, by (3) we may choose a subsequence which converges (uniformly on compacts) to some non-trivial function $f\in B_\sigma.$ Clearly, $f=0$ on $\La$. We are back to the previous  case.

\medskip
\noindent{\bf 3}.  Assume $\|f_n\cdot 1_N\|_\infty\to 0, n\to\infty,$ for every $N>0$. Let $x_n$ be a point such that $|f(x_n)|>1/2.$ Passing to a subsequence and using (\ref{ber}), we may assume that $f_n$ is so small on $[-2n,2n]$ that the classical Two-constants theorem (see \cite{n}, ch. 3) implies
$$
|f_n(z)|<\frac{1}{n^3}, \ \ |z|<n.
$$
From this and (\ref{ber}) it follows  that the series $$f(z):=\sum_{n\in\N}nf_n(z)$$uniformly converges on compacts and that $f$ belongs to $E_\sigma$. Clearly, $f$ satisfies $|f(x_n)|>n/2,$ so that it is not bounded on $\R$. On the other hand, by (\ref{sist}), we see that $f$ is bounded on $\La.$

\section{Proof of Lemma 2}
Here we present a sketch of a proof of Lemma 2.

Given $a\ne0$ and $n\in\N$, set
$$
P_{a,n}(z):=\left(1-\frac{z^2}{a^2}\right)^{n-1}.
$$
We will use

\begin{lemma}
Assume $f\in B_\sigma, \|f\|_\infty\leq 1.$  Then for every $\varepsilon>0$ and $x_0>0$ there exist $a>1$ and $n\in\N$ such that
$$
\|P_{a,n}f\|_\infty \leq 1, \ \|P_{a,n+1}f\|_\infty>1, \ |P_{a,n}(x)|>1-\varepsilon, \ \ |x|\leq x_0.
$$
\end{lemma}

 The proof follows relatively easily from  the following  fact: If $f\in B_\sigma$, then for every $\varepsilon>0$ the set
$$
\{x\in\R: |f(x)|>e^{-\varepsilon |x|}\}
$$is unbounded
 (see, for example, \cite{levin}, Lec. 16, Theorem 2). Note that  the latter estimate is a consequence of the fundamental fact that for every $f\in B_\sigma$ the so-called logarithmic integral of $f$ converges, see \cite{levin}, Lec.~16.

\medskip\noindent{\bf  Proof of Lemma 2}.
Since $\La$ is not a uniqueness set for $B_\sigma,$ there is a function $f\in  B_\sigma$ which vanishes on $\La$. We may assume that $\|f\|_\infty\leq 1$.

Let us consider two cases.

\medskip\noindent{\bf 1}. Assume there is an integer $n$ such that
$$
\limsup_{|x|\to\infty}|x|^n|f(x)|=\infty.
$$
Then the function $g(z):=z^{n}f(z)\in E_\sigma$, vanishes on $\La$ and is unbounded on $\R$.

\medskip\noindent{\bf 2}.  Assume
$$
\limsup_{|x|\to\infty}|x|^n|f(x)|=0, \ \ \mbox{ for every } n>0.
$$

By Lemma 3, there exist  $a_1>1$ and $n_1$ such that the function  $$f_1(z):=P_{a_1,n_1}(z)f(z)$$ satisfies
  $\|f_1\|_\infty\leq 1$ and $\|(1-x^2/a_1^2)f_1(x)\|_\infty>1$. The latter inequality shows that there is a point $x_1$ such that
\begin{equation}\label{t}
|f_1(x_1)|\geq \frac{1}{x_1^2}.
\end{equation}Clearly, $f_1\in B_\sigma.$

Similarly, we construct a function $$f_2(z):=P_{a_1,n_1}(z)f_1(z)\in B_\sigma,$$ such that
$ \|f_2\|_\infty\leq 1$ and there is a point $x_2$ such that
$$
|f_2(x_2)|\geq \frac{1}{x_2^2}.
$$In addition, we may assume that $|P_{a_2,n_2}(x)|>1/2$  for $ |x|\leq x_1,$ so that we have
$$
|f_2(x_1)|> \frac{1}{2x_1^2}.
$$

One may repeat  this procedure $n$ times to construct a sequence of functions $f_n\in B_\sigma$ satisfying $\|f_n\|_\infty\leq 1$. Moreover, on each step we may assume that $|P_{a_k,n_k}(x)|$ is so close to one for $|x|\leq x_{k-1}$, that we have
\begin{equation}\label{ttt}
|f_n(x_k)|> \frac{1}{2x_k^2}, \ \ 1\leq k\leq n.
\end{equation}
By (2), the
 family $f_n\in B_\sigma$ is normal. Also, every function $f_n$ vanishes on $\La$. Hence, there is a subsequence $n=n_j\to\infty$ such that $f_n$ converge uniformly on  compacts to some function $\varphi\in B_\sigma$. Clearly, the function $g(z):=z^3\varphi(z)$ belongs to $E_\sigma$, and it follows from (\ref{ttt}) that it is unbounded on $\R$.

%\medskip Observe, that a similar approach can be used to prove Proposition~3.

\medskip\noindent{\bf Acknowledgement}. The authors are grateful to Misha Sodin for the constructive comments.

\medskip

 Stavanger University, 4036 Stavanger,  Norway.

 \medskip

 E-mail addresses:

   natalia.blank@uis.no

 alexander.ulanovskii@uis.no
\end{large}
\end{document}